\documentclass{article}%
\usepackage{amssymb}
\usepackage{graphicx}
\usepackage{amsmath}
\usepackage{amsfonts}%
\providecommand{\U}[1]{\protect\rule{.1in}{.1in}}
\newtheorem{theorem}{Theorem}

\newtheorem{corollary}[theorem]{Corollary}

\newtheorem{lemma}[theorem]{Lemma}

\newtheorem{proposition}[theorem]{Proposition}
\newtheorem{remark}[theorem]{Remark}

\begin{document}
\textbf{To appear in Probability Theory and Related Fields}

\begin{center}
{\LARGE Central limit theorem for triangular arrays of Non-Homogeneous Markov
chains\ }

\bigskip Magda Peligrad\footnote{Supported in part by a Charles Phelps Taft
Memorial Fund grant, NSF\ DMS-0830579 and NSA\ grants H98230-09-1-0005 and
H98230-11-1-0135.}
\end{center}

\bigskip\textbf{ }

Department of Mathematical Sciences, University of Cincinnati, PO Box 210025,
Cincinnati, Oh 45221-0025, USA. Email: peligrm@ucmail.uc.edu \vskip10pt

\bigskip

\textbf{Abbreviated Title: CLT for non-homogeneous Markov chains}

\bigskip



{}

\begin{center}
{\Large Abstract}
\end{center}

In this paper we obtain the central limit theorem for triangular arrays of
non-homogeneous Markov chains under a condition imposed to the maximal
coefficient of correlation. The proofs are based on martingale techniques and
a sharp lower bound estimate for the variance of partial sums. The results
complement an important central limit theorem of Dobrushin based on the
contraction coefficient.

\section{Introduction and notations}

\qquad More than fifty years ago Dobrushin \cite{Do2} proved a definitive
central limit theorem for non-homogeneous Markov chains. His work is based on
the coefficient of ergodicity which is defined by using the contraction
coefficient, specifically for uniformly bounded functions. In a recent paper,
Sethuraman and Varadhan \cite{SV} give a new and elegant proof of Dobrushin's
result and provide a survey of the literature that was generated by it. In
this paper we address a similar problem for Markov chains by using the maximal
coefficient of correlation, instead of the contraction coefficient. This
coefficient is more general and the results are applicable to a larger class
of Markov processes. The problem is challenging, since the maximal coefficient
of correlation is defined for functions that are square integrable only and
many new tools have to be developed.

\bigskip

Let $(\Omega,\mathcal{K},\mathbf{P})$ be a probability space and let
$\mathcal{A},\mathcal{B}$ be two sub $\sigma$-algebras of $\mathcal{K}$.
Define the maximal coefficient of correlation
\[
\rho(\mathcal{A},\mathcal{B})=\sup_{f\in\mathbf{L}_{2}(\mathcal{A}%
),g\in\mathbf{L}_{2}(\mathcal{B})}|corr\,(f,g)|\text{ ,}%
\]
where $\mathbf{L}_{2}(\mathcal{A})$ is the space of random variables that are
$\mathcal{A}$ measurable and square integrable. For a vector of random
variables, $(Y_{k})_{1\leq k\leq n}$ we define
\begin{equation}
\rho_{k}=\max_{1\leq s,s+k\leq n}\rho(\sigma(Y_{i},i\leq s),\sigma(Y_{j},j\geq
s+k))\text{ .} \label{defrogen}%
\end{equation}
For a nonhomogeneous Markov chain of length $n$, $(\xi_{i})_{1\leq i\leq n}$,
it turns out that the computation of $\rho_{k}$ simplifies (see for instance
Theorem 7.2 (c) in \cite{br}). For this case,
\[
\rho_{k}=\max_{1\leq s,s+k\leq n}\rho(\sigma(\xi_{s}),\sigma(\xi_{s+k}))\text{
.}%
\]
Moreover (see Theorem 7.4 (a) in \cite{br}), for all $1\leq k\leq n-1,$
\[
\rho_{k}\leq\rho_{1}^{k}\text{ .}%
\]
In terms of the conditional expectation (see chapter 7 in \cite{Ros} or
Theorem 4.4 (b3) in \cite{br}) an alternative definition of $\rho_{1}$ is
\begin{equation}
\rho_{1}=\max_{2\leq i\leq n}\sup_{g}{\LARGE \{}\frac{||\mathbf{E}(g(\xi
_{i})|\xi_{i-1})||_{2}}{||g(\xi_{i})||_{2}}\text{ ; }||g(\xi_{i})||_{2}%
<\infty\text{ and }\mathbf{E}g(\xi_{i})=0{\LARGE \}}\text{ ,} \label{def-ro}%
\end{equation}
where we used the notation $||X||_{p}=(\mathbf{E}|X|^{p})^{1/p}$, for $p>1$.

For a stationary Markov chain defined on $(\Omega,\mathcal{K},P)$ with values
in $(\mathcal{X},\mathcal{B}(\mathcal{X}))$ with invariant measure $\pi$ and
transition probability $Q(x,A)=P(\xi_{1}\in A|\xi_{0}=x)$, define the operator
$Q$ acting on $\mathbf{L}_{2}(\mathcal{X},\mathcal{B}(\mathcal{X}),\pi)$ via%
\begin{equation}
(Qu)(x)=\int_{\mathcal{X}}u(y)Q(x,dy)\text{ .} \label{opQ}%
\end{equation}
Denote $\mathbf{L}_{2}^{0}(\pi)=\{g\in\mathbf{L}_{2}(\mathcal{X}%
,\mathcal{B}(\mathcal{X}),\pi)$ with $\int gd\pi=0\}.$ With these notations,
the coefficient $\rho_{1}$ is simply the norm operator of $Q:\mathbf{L}%
_{2}^{0}(\pi)\rightarrow\mathbf{L}_{2}^{0}(\pi)$,
\begin{equation}
\rho_{1}=||Q||_{\mathbf{L}_{2}^{0}(\pi)}=\sup_{g\in\mathbf{L}_{2}^{0}(\pi
)}\frac{||Q(g)||_{2}}{||g||_{2}}\text{.} \label{ro-sta}%
\end{equation}

Conditions imposed to the maximal coefficient of correlation make possible to
study the asymptotic behavior of many dependent structures including classes
of Markov chains and Gaussian sequences. This coefficient was used by
Kolmogorov and Rozanov \cite{KR} and further studied by Rosenblatt \cite{Ros},
Ibragimov \cite{Ib}, Shao \cite{Shao} among many others. An introduction to
this topic, mostly in the stationary setting, can be found in the Chapters 7,
9 and 11 in Bradley \cite{br}. Application to the central limit theorem (CLT)
for various stationary Markov chains with $\rho_{1}<1$ are surveyed in Jones
\cite{Jo}. In the nonstationary setting and general triangular arrays a
central limit theorem was obtained by Utev \cite{Ut}, assuming a lower bound
on the variance of partial sums and $\rho-$mixing coefficients converging to
$0$ uniformly at a logarithmic rate.

In this paper we are concerned with the central limit theorem for a triangular
array of Markov chains. Let $(\xi_{n,i})_{1\leq i\leq n}$ be an array of
non-homogeneous Markov chains defined on a probability space $(\Omega
,\mathcal{K},P)$ with values in $(\mathcal{X},\mathcal{B}(\mathcal{X})).$ In
addition, let $(f_{n,i})_{1\leq i\leq n}$ be real valued functions on
$\mathcal{X}$. Define,
\begin{equation}
X_{n,i}=f_{n,i}(\xi_{n,i})\text{ and }S_{n}=\sum_{i=1}^{n}X_{n,i}\text{ .}
\label{def X}%
\end{equation}
Everywhere in the paper we shall assume%
\[
\mathbf{E}X_{n,i}=0,\text{ }\mathbf{E}X_{n,i}^{2}<\infty
\]
and denote by%
\begin{equation}
\sigma_{n}^{2}=var\text{ }S_{n}\text{ and \ }b_{n}^{2}=\sum_{i=1}^{n}var\text{
}X_{n,i}\text{ ,} \label{def b}%
\end{equation}
where $var$ $X=\mathbf{E}(X-\mathbf{E}X)^{2}.$ In this context, the mixing
coefficients need an additional index to indicate the row. We shall write now
$\rho_{n,1}$ instead of $\rho_{1}$ to specify that this coefficient is
computed for $(\xi_{n,i})_{1\leq i\leq n}.$

We shall establish in Proposition \ref{lowerB2} that for all $n\geq1,$ the
quantities $\sigma_{n}^{2}$ and $b_{n}^{2}$ are related by the following
inequality
\[
\frac{1-\rho_{n,1}}{1+\rho_{n,1}}b_{n}^{2}\leq\sigma_{n}^{2}\leq\frac
{1+\rho_{n,1}}{1-\rho_{n,1}}b_{n}^{2}\text{ ,}%
\]
provided $\rho_{n,1}<1.$

We shall further discuss the rate at which the maximal coefficients of
correlation, $\rho_{n,1}$, are allowed to converge to $1$, for the validity of
the central limit theorem for $S_{n}/\sigma_{n}$.

The results are formulated both in terms of bounded random variables and also
in an integral form similar to the Lindeberg condition. When applied to a
triangular array of uniformly bounded random variables, with the variance of
individual summands uniformly bounded below, our sufficient condition for the
central limit theorem is implied, for instance, by
\begin{equation}
(1-\rho_{n,1})^{3}n(\ln n)^{-2}\rightarrow\infty\text{ as }n\rightarrow
\infty\text{ .} \label{log1}%
\end{equation}
We can see from this result that we obtain the central limit theorem not only
for the situation when $\rho_{n,1}\leq r<1.$ We can let $\rho_{n,1}$
approaches $1,$ but not too fast, at a rate that will be specified.

The proof of this result and the other results of this type are based on the
following tools we develop in this paper that have interest in themselves:

(1) General sufficient conditions for the CLT for triangular arrays, based on
a familiar projective martingale representation.

(2) Sharp lower and upper bounds for the variance of sums of variables
connected in a Markov chain.

(3) Moment and exponential inequalities for certain partial sums.

Our theorems are related to Dobrushin's result which uses as a measure of
dependence the contraction coefficient
\[
\delta(Q)=\sup_{u\in\mathcal{U}}\sup_{x_{1},x_{2}}|(Qu)(x_{1})-(Qu)(x_{2}%
)|\text{ ,}%
\]
where $\mathcal{U}=\{u$, $\sup_{x_{1},x_{2}}|u(x_{1})-u(x_{2})|\leq1\}$ and
the operator $Q$ is defined by (\ref{opQ}), on the space of bounded measurable
functions. For a triangular array of Markov chains $(\xi_{n,i})_{1\leq i\leq
n}$ with transition probabilities $Q_{n,i}(x,A)=P(\xi_{n,i+1}\in A|\xi
_{n,i}=x)$,
\begin{equation}
\delta_{n,1}=\sup_{1\leq i\leq n-1}\delta(Q_{n,i})\text{ and }\delta_{n,k}%
\leq\delta_{n,1}^{k}\text{.} \label{defdelta}%
\end{equation}
By Lemma 4.1 in Sethuraman and Varadhan \cite{SV} we have%
\[
\rho_{n,1}\leq\delta_{n,1}^{1/2}\text{ .}%
\]
Dobrushin \cite{Do2} showed that for a triangular array of uniformly bounded
random variables, with the variance of individual summands uniformly bounded
below, a sufficient condition for the CLT is
\begin{equation}
(1-\delta_{n,1})^{3}n\rightarrow\infty\text{ as }n\rightarrow\infty\text{ .}
\label{contract}%
\end{equation}
Moreover, he analyzed an example attributed to Bernstein, of family of Markov
chains satisfying $1-\delta_{n,1}=n^{1/3}$ and such that the CLT fails. The
initial proof of this result is very long. A simplified proof of a further
reaching result can be found in Sethuraman and Varadhan \cite{SV} (see also
Theorem \ref{Dobr}). For the properties of this contraction coefficient we
refer to Iosifescu and Theodorescu \cite{Io2}, sections 1.1. and 1.2.

There are plenty of examples for which $\delta_{n,1}=1,$ but $\rho_{n,1}<1,$
so our results have a larger sphere of applicability than the results based on
$\delta_{n,1}$.

For instance, for a row-wise stationary array of Markov chains with joint
distribution of $(\xi_{n,1},\xi_{n,2})$ bivariate normal, the maximal
coefficient of correlation is very simple, namely $\rho_{n,1}=|corr(\xi
_{n,1},\xi_{n,2})|$, while (if i.i.d. rows are excluded) $\delta_{n,1}=1$ for
all $n$ (for a convenient reference to this fact see \cite{br}, Theorems 9.1
and 9.7). For functions of variables in this array our theorems are
applicable, and the conditions are very easy to verify.

Moreover, even for the situation when $\delta_{n,1}<1,$ it is possible that
$\rho_{n,1}\rightarrow1$ and $\delta_{n,1}\rightarrow1$ at different rates,
such that, for instance, (\ref{log1}) holds but (\ref{contract}) does not.
Such an example can be easily constructed by using a recent result by Bradley
\cite{BR}, who showed that for any $0<a<b<1$ there is a stationary Markov
chain for which $\rho_{1}=a$ and $\delta_{1}=b.$

Our results will be useful for treating families of various Markov processes
that are considered in applications. For example, Liu et all \cite{Liu} have
shown that if the operator induced by a Gibbs sampler satisfies a
Hilbert-Schmidt condition then $\rho_{1}<1$.

Another class of examples is provided by an array of stationary reversible
Markov chains that are geometrically ergodic. A stationary Markov chain is
called geometrically ergodic if there is $0<t<1$ and a nonnegative function
$M(x)$ such that $||Q^{n}(x,\cdot)-\pi(\cdot)||\leq M(x)t^{n}$. A stationary
Markov chain that is geometrically ergodic and reversible satisfies $\rho
_{1}<1$. (Roberts and Rosenthal, \cite{RR}). A particular example of this kind
is the popular Random Walk MHG Algorithms\textbf{ }which are reversible by
construction. In Mengersen and Tweedie \cite{Ma} it was shown that the random
walk samplers cannot be uniformly ergodic (so $\delta_{1}=1$) but they do
establish that a random walk MHG algorithm can be geometrically ergodic in
some situations, therefore they have $\rho_{1}<1$. This work was extended in
Roberts and Tweedie \cite{RT}.

Our paper is organized as follows: In Section 2 we state the main results. In
order to prove them, in Section 3 we develop sufficient conditions for the CLT
for triangular arrays of random variables based on martingale representations.
Section 4 is concerned with bounds for the variance of partial sums of a
Markov chain as a function of the $\rho_{1}$ coefficient. The proofs of the
main results are the subject of Section 5. Some technical lemmas involving
higher moments for sums and exponential bounds are postponed to the Appendix.

The convergence in probability will be denoted $\rightarrow^{P}$ and
$\rightarrow^{D}$ denotes convergence in distribution.

\section{Results}

Our first theorem applies to triangular arrays of functions of Markov chains
consisting of bounded centered variables. To describe our results it is
convenient to introduce the related coefficient%
\begin{equation}
\lambda_{n}=1-\rho_{n,1}=\min_{1\leq s\leq n-1}[1-\rho(\sigma(\xi
_{n,s}),\sigma(\xi_{n,s+1}))]\text{ .} \label{def-lam}%
\end{equation}
Clearly $0\leq\lambda_{n}\leq1$, and $\lambda_{n}$ is a coefficient of
independence for $(\xi_{n,i})_{1\leq i\leq n}$ with $n$ fixed. Notice that
$\lambda_{n}=1$ if and only if the vector $(X_{n,i})_{1\leq i\leq n}$ is
independent. Everywhere in this section we shall consider the nondegenerate
case i.e. $\lambda_{n}b_{n}>0$ for all $n\geq1.$ By Proposition \ref{lowerB2}
this condition is equivalent to $\lambda_{n}\sigma_{n}>0$ for all $n\geq1.$

\begin{theorem}
\label{CLTMarkov}Suppose that $(X_{n,i})_{1\leq i\leq n}$ is defined by
(\ref{def X}) and for some finite positive constants $C_{n}$ we have
\begin{equation}
\max_{1\leq i\leq n}|X_{n,i}|\leq C_{n}\text{ a.s.} \label{bound}%
\end{equation}
and
\begin{equation}
\frac{C_{n}(1+|\ln(\lambda_{n})|)}{\lambda_{n}\sigma_{n}}\rightarrow0\text{
\ as }n\rightarrow\infty\text{ }. \label{Dob}%
\end{equation}
Then%
\begin{equation}
\frac{\sum_{i=1}^{n}X_{n,i}}{\sigma_{n}}\overset{\mathcal{D}}{\rightarrow
}N(0,1)\text{ as \ }n\rightarrow\infty\text{ }. \label{CLTD}%
\end{equation}

\end{theorem}

We state a corollary which combines Theorem \ref{CLTMarkov} with the bound of
the variance given in Proposition \ref{lowerB2} in Section 4.

\begin{corollary}
Suppose that $(X_{n,i})_{1\leq i\leq n}$ is defined by (\ref{def X}). Assume
that (\ref{bound}) holds and
\[
\frac{C_{n}(1+|\ln(\lambda_{n})|)}{\lambda_{n}^{3/2}b_{n}}\rightarrow0\text{
as }n\rightarrow\infty\text{ .}%
\]
Then the CLT (\ref{CLTD}) holds.
\end{corollary}

Next, we give a corollary that can be applied to an array of uniformly bounded
random variables.

\begin{corollary}
Suppose that $(X_{n,i})_{1\leq i\leq n}$ is defined by (\ref{def X}) and
assume that there are two positive constants $C$ and $c$ such that
$\max_{1\leq i\leq n}|X_{n,i}|\leq C$ a.s. and also $var$ $X_{n,i}\geq c>0$
for all $n\geq1$ and $1\leq i\leq n.$ Then CLT (\ref{CLTD}) holds provided
\begin{equation}
\lambda_{n}^{3}n(1+|\ln(\lambda_{n})|)^{-2}\rightarrow\infty. \label{log2}%
\end{equation}

\end{corollary}

Notice that (\ref{log1}) implies (\ref{log2}).

\bigskip

We shall also prove an integral form of Theorem \ref{CLTMarkov}.

\begin{corollary}
\label{Lindeberg}Suppose that $(X_{n,i})_{1\leq i\leq n}$ is defined by
(\ref{def X}) and for every $\varepsilon>0$
\begin{equation}
\frac{1}{\lambda_{n}\sigma_{n}^{2}}\sum_{i=1}^{n}\mathbf{E}X_{n,i}%
^{2}I(|X_{n,i}|>\varepsilon h(\lambda_{n})\sigma_{n})\rightarrow0\text{ as
}n\rightarrow\infty\label{lindeberg}%
\end{equation}
where $h(\lambda_{n})=\lambda_{n}(1+|\ln(\lambda_{n})|)^{-1}.$ Then the CLT
(\ref{CLTD}) holds.
\end{corollary}

We now point out some immediate consequences of Corollary \ref{Lindeberg}.

\begin{remark}
\label{remark5}Let us notice that by using the bounds on the variance given in
Proposition \ref{lowerB2}, condition (\ref{lindeberg}) is implied by%
\begin{equation}
\frac{1}{\lambda_{n}^{2}b_{n}^{2}}\sum_{i=1}^{n}\mathbf{E}X_{n,i}%
^{2}I(|X_{n,i}|>\varepsilon h^{\prime}(\lambda_{n})b_{n})\rightarrow0\text{ as
}n\rightarrow\infty\label{Lindebergb}%
\end{equation}
where $h^{\prime}(\lambda_{n})=\lambda_{n}^{3/2}(1+|\ln(\lambda_{n})|)^{-1}.$
\end{remark}

The next remark applies to triangular arrays of Markov chains with uniformly
bounded $\rho_{n,1}-$mixing coefficients.

\begin{remark}
Suppose that $(X_{n,i})_{1\leq i\leq n}$ is defined by (\ref{def X}) and there
is a positive number $\rho$ such that $\sup_{n}\rho_{n,1}\leq\rho<1.$ Then,
the CLT (\ref{CLTD}) holds provided that for every $\varepsilon>0$
\[
\frac{1}{b_{n}^{2}}\sum_{i=1}^{n}\mathbf{E}X_{n,i}^{2}I(|X_{n,i}|>\varepsilon
b_{n})\rightarrow0\text{ .}%
\]

\end{remark}

For arrays of Markov chains that are row-wise strictly stationary with the
same invariant distribution, Corollary \ref{Lindeberg} (via Remark
\ref{remark5}) has a simple form. Examples of this type are arrays of Markov
chains generated by parametric copulas.

\begin{remark}
Suppose for each $n$, $(\xi_{n,i})_{1\leq i\leq n}$ is a stationary Markov
chain with the same invariant distribution $\pi$ and transition operator
$Q_{n}$. Let $f\in\mathbf{L}_{2}^{0}(\pi)$ and define $X_{n,k}=f(\xi_{n,k})$
and $\lambda_{n}=1-||Q_{n}||_{\mathbf{L}_{2}^{0}(\pi)}$. Then, the CLT
(\ref{CLTD}) holds provided that for every $\varepsilon>0$
\[
\frac{1}{\lambda_{n}^{2}}\int f^{2}(x)I(|f(x)|>\varepsilon\sqrt{n}h^{\prime
}(\lambda_{n}))d\pi\rightarrow0\text{ as }n\rightarrow\infty\text{ .}%
\]

\end{remark}

Finally, we mention the extension of Dobrushin's CLT obtained by Sethuraman
and Varadhan \cite{SV}, by using the coefficient $\delta_{n,1}$:

\begin{theorem}
\label{Dobr}Suppose that $(X_{n,i})_{1\leq i\leq n}$ is defined by
(\ref{def X}) and relation (\ref{bound}) holds. Denote $\alpha_{n}%
=1-\delta_{n,1}.$ If
\begin{equation}
\frac{C_{n}^{2}}{\alpha_{n}^{3}b_{n}^{2}}\rightarrow0\text{ as }%
n\rightarrow\infty\text{ ,} \label{CD}%
\end{equation}
then, the CLT (\ref{CLTD}) holds. When $|X_{n,i}|\leq C<\infty$ a.s and $var$
$X_{n,i}\geq c>0$, for all $1\leq i\leq n$ and $n\geq1$, then the CLT
(\ref{CLTD}) holds provided
\[
\alpha_{n}^{3}n\rightarrow\infty\text{ as }n\rightarrow\infty\text{ .}%
\]

\end{theorem}

\section{Central Limit Theorem for triangular arrays}

The following theorem is a variant of Theorem 3.2 in Hall and Heyde \cite{HH}
(see also G\"{a}nssler and H\"{a}usler \cite{GH}).

\begin{theorem}
\label{martCLT}Assume $(D_{n,i})_{1\leq i\leq n}$ is an array of square
integrable martingale differences adapted to an array $(\mathcal{F}%
_{n,i})_{1\leq i\leq n}$ of nested sigma fields. Suppose
\begin{equation}
\mathbf{E}(\max_{1\leq j\leq n}|D_{n,j}|)\rightarrow0\text{ as }%
n\rightarrow\infty\label{negl1}%
\end{equation}
and%
\[
\sum_{j=1}^{n}D_{n,j}^{2}\rightarrow^{P}1\text{ }\ \text{as }n\rightarrow
\infty\text{ .}%
\]
Then $S_{n}=\sum_{j=1}^{n}D_{n,j}$ converges in distribution to a standard
normal variable.
\end{theorem}

In this section we shall assume the following general setting:

\bigskip

$(C1)$ Assume $(X_{n,j})_{1\leq j\leq n}$ is an array of centered random
variables that are square integrable and adapted to an array of sigma fields
$(\mathcal{F}_{n,j})_{1\leq j\leq n},$ with $\mathcal{F}_{n,j}\subset
\mathcal{F}_{n,j+1}$ for all $1\leq j\leq n-1$. Extend the array with
$X_{n,0}=0$ and $\mathcal{F}_{n,0}=\{\varnothing,\Omega\}$ for all $n.$

With this notation, as an immediate consequence of Theorem \ref{martCLT}, we formulate:

\begin{corollary}
\label{proj}Assume $(C1)$ and let $\mathbf{E}S_{n}^{2}=1.$ Define the
projector operator$\ $%
\[
\mathbf{P}_{n,j}Y=\mathbf{E}(Y|\mathcal{F}_{n,j})-\mathbf{E}(Y|\mathcal{F}%
_{n,j-1})\text{ .}%
\]
Assume%
\begin{equation}
\max_{1\leq j\leq n}|\mathbf{P}_{n,j}S_{n}|\rightarrow^{P}0\text{ as
}n\rightarrow\infty\text{ .} \label{NEG}%
\end{equation}
and%
\begin{equation}
\sum_{j=1}^{n\ }(\mathbf{P}_{n,j}S_{n})^{2}\rightarrow^{P}1\text{ as
}n\rightarrow\infty\text{ .} \label{Conv}%
\end{equation}
Then the CLT in (\ref{CLTD}) holds.
\end{corollary}

\textbf{Proof}. Because we assume $\mathcal{F}_{n,0}=\{\varnothing,\Omega\},$
we can express $S_{n}$ in terms of projections
\begin{equation}
S_{n}=\sum_{j=1}^{n}X_{n,j}=\sum_{j=1}^{n}\mathbf{E}(S_{n}|\mathcal{F}%
_{n,j})-\mathbf{E}(S_{n}|\mathcal{F}_{n,j-1})=\sum_{j=1}^{n}\mathbf{P}%
_{n,j}S_{n}\text{ .} \label{martdec2}%
\end{equation}
Notice that we have written $S_{n}$ as a sum of martingale differences
\begin{equation}
d_{n,j}=\mathbf{E}(S_{n}|\mathcal{F}_{n,j})-\mathbf{E}(S_{n}|\mathcal{F}%
_{n,j-1})=\mathbf{P}_{n,j}S_{n}=\mathbf{P}_{n,j}(S_{n}-S_{j-1}) \label{defD}%
\end{equation}
and we apply Theorem \ref{martCLT}. Since $\sum_{j=1}^{n\ }\mathbf{E(P}%
_{n,j}S_{n})^{2}=\mathbf{E}S_{n}^{2}=1,$ it follows that $\max_{1\leq j\leq
n}|\mathbf{P}_{n,j}S_{n}|$ is uniformly integrable in $\mathbf{L}_{1}$ and
then (\ref{NEG}) implies (\ref{negl1}). $\Diamond$

\bigskip

Analyzing the conditions of Corollary \ref{proj} is leading us to the
following useful theorem. For $0\leq j\leq n$ denote%
\begin{equation}
A_{n,j}=\mathbf{E}(S_{n}-S_{n,j}|\mathcal{F}_{n,j})\text{ ,} \label{Def A}%
\end{equation}
where $S_{n,j}=\sum_{i=1}^{j}X_{n,i}$

\begin{theorem}
\label{CLT} Assume $(C1)$ and $\mathbf{E}S_{n}^{2}=1.$ Also assume that
\begin{equation}
\max_{1\leq j\leq n}(|X_{n,j}|+|A_{n,j}|{\large )}\rightarrow^{P}0\text{ as
}n\rightarrow\infty\text{ } \label{NEGT}%
\end{equation}
and%
\begin{equation}
\sum_{j=1}^{n}{\large (}X_{n,j}^{2}+2X_{n,j}A_{n,j}{\large )}\rightarrow
^{P}1\text{ as }n\rightarrow\infty\text{ .} \label{convt}%
\end{equation}
Then $S_{n}$ converges in distribution to $N(0,1)$.
\end{theorem}

\textbf{Proof.} For simplicity we drop the index $n$ in the notation, so
$X_{j}=X_{n,j}$, $A_{j}=A_{n,j},$ $d_{j}=d_{n,j}$.

Condition (\ref{NEG}) follows from condition (\ref{NEGT}) since, by
definitions (\ref{defD}) and (\ref{Def A}),
\[
|d_{j}|\leq|X_{j}|+|A_{j}|+|A_{j-1}|\text{ a.s.}%
\]
To verify condition (\ref{Conv}), for $1\leq j\leq n$ let us compute
\begin{align*}
d_{j}^{2}  &  =(X_{j}+A_{j}-A_{j-1})^{2}\\
&  =(X_{j}^{2}+2X_{j}A_{j})+A_{j}^{2}-A_{j-1}^{2}+2(A_{j-1}-X_{j}%
-A_{j})A_{j-1}\text{ ,}%
\end{align*}
whence, by definition (\ref{defD}) and the fact that $A_{0}=A_{n}=0,$ we
obtain
\[
\sum_{j=1}^{n}d_{j}^{2}=\sum_{j=1}^{n}(X_{j}^{2}+2X_{j}A_{j})-2\sum_{j=1}%
^{n}d_{j}A_{j-1}\text{ .}%
\]
Now, by condition (\ref{convt}), the first term in the right hand side is
converging in probability to $1$. For the martingale transform $\sum_{j=1}%
^{n}d_{j}A_{j-1}$ we use a truncation argument. Let $\varepsilon>0$ and denote
$A_{j\ }^{\varepsilon}=A_{j}I(|A_{j}|\leq\varepsilon).$ For any $a>0,$
\begin{gather*}
\mathbf{P}(|\sum_{j=1}^{n}d_{j}A_{j-1}|>a)\leq\mathbf{P}(\max_{1\leq j\leq
n}|A_{j}|>\varepsilon)+\mathbf{P}(|\sum_{j=1}^{n}d_{j}A_{j-1}^{\varepsilon
}|>a)\\
\leq\mathbf{P}(\max_{1\leq j\leq n}|A_{j}|>\varepsilon)+\varepsilon
^{2}\mathbf{E}(\sum_{j=1}^{n}d_{j})^{2}/a^{2}=\mathbf{P}(\max_{1\leq j\leq
n}|A_{j}|>\varepsilon)+\varepsilon^{2}/a^{2}.
\end{gather*}
where on the last line we used the fact that $\mathbf{E}(\sum_{j=1}^{n}%
d_{j})^{2}=1.$ Then, we take into account that $\max_{1\leq j\leq n}|A_{j}|$
is negligible in probability by (\ref{NEGT}) and we conclude the convergence
to $0$ by letting $n\rightarrow\infty$ followed by $\varepsilon\rightarrow0.$
It follows that
\[
\sum_{j=1}^{n}d_{j}^{2}\rightarrow^{P}1\text{ as }n\rightarrow\infty
\]
and the CLT holds by Corollary \ref{proj}. $\Diamond$

\bigskip

Theorem \ref{CLT} has the following simple Corollary that will be used in our proofs:

\begin{proposition}
\label{markov} Assume $(C1)$ and the variables have finite moments of order
$4$. Moreover assume the following conditions hold:%
\begin{equation}
\frac{1}{\sigma_{n}^{4}}\sum_{j=1}^{n}\mathbf{E}X_{n,j}^{4}\rightarrow0\text{
as }n\rightarrow\infty\text{ .} \label{C1}%
\end{equation}
For every $\varepsilon>0$%
\begin{equation}
\mathbf{P}(\max_{1\leq j\leq n}|A_{j}|>\varepsilon\sigma_{n})\rightarrow
0\text{ as }n\rightarrow\infty\label{C2}%
\end{equation}
and%
\begin{equation}
\frac{1}{\sigma_{n}^{4}}var\sum_{j=1}^{n}(X_{n,j}^{2}+2X_{n,j}A_{n,j}%
)\rightarrow0\text{ as }n\rightarrow\infty\text{ .} \label{C3}%
\end{equation}
Then $\sigma_{n}^{-1}S_{n}$ converges in distribution to $N(0,1)$ .
\end{proposition}

\textbf{Proof}. Conditions (\ref{C1}) and (\ref{C2}) easily imply
(\ref{NEGT}). Then, condition (\ref{C3}) implies (\ref{convt}) by taking into
account that
\[
\frac{1}{\sigma_{n}^{2}}\mathbf{E}\sum_{j=1}^{n}(X_{n,j}^{2}+2X_{n,j}%
A_{n,j})=1\text{ .}%
\]
$\Diamond$

\section{Bounds for the variance of partial sums of Markov chains}

In this section we establish sharp upper and lower bounds for the variance of
partial sums of a Markov chain as a function of the maximal coefficient of
correlation defined in (\ref{def-ro}).

\begin{proposition}
\label{lowerB2}Let $(X_{1},X_{2},...,X_{n})$ be a vector of square integrable
centered random variables that are functions of a Markov process $(\xi
_{i})_{1\leq i\leq n}$ i.e. $X_{k}=f_{k}(\xi_{k})$. Denote by $S_{n}%
=\sum_{i=1}^{n}X_{i}$ and $\mathcal{F}_{j}=\sigma(\xi_{i},$ $i\leq j).$ We set
$X_{0}=0$ and $\mathcal{F}_{0}=\{0,\Omega\}.$ If $\rho_{1}<1,$ then%
\[
\frac{1-\rho_{1}}{1+\rho_{1}}\sum_{i=1}^{n}\mathbf{E}X_{i}^{2}\leq
\mathbf{E}S_{n}^{2}\leq\frac{1+\rho_{1}}{1-\rho_{1}}\sum_{i=1}^{n}%
\mathbf{E}X_{i}^{2}\text{ .}%
\]

\end{proposition}

\textbf{Proof.} We prove first the lower bound. For this proof we recall the
notation (\ref{def b}) and (\ref{Def A}), which in this case is $A_{j}%
=A_{n,j}=\mathbf{E}(S_{n}-S_{j}|\mathcal{\xi}_{j})$. We start as before from
the martingale decomposition
\[
S_{n}=\sum_{j=1}^{n}\mathbf{E}(S_{n}|\mathcal{F}_{j})-\mathbf{E}%
(S_{n}|\mathcal{F}_{j-1})=\sum_{j=1}^{n\ }\mathbf{P}_{j}(S_{n}-S_{j-1})\text{
.}%
\]
By the orthogonality of the martingale differences%
\begin{equation}
\sigma_{n}^{2}=\sum_{j=1}^{n\ }\mathbf{E}[\mathbf{P}_{j}(S_{n}-S_{j-1}%
)]^{2}\text{ .} \label{sigma}%
\end{equation}
Notice that by taking into account the Markov property and simple algebra
{\LARGE \ }%
\begin{gather*}
\mathbf{E}[\mathbf{P}_{j}(S_{n}-S_{j-1})]^{2}=\mathbf{E}(X_{j}+A_{j}%
)^{2}+\mathbf{E}(A_{j-1})^{2}\\
-2\mathbf{E}(X_{j}+A_{j})A_{j-1}\text{ .}%
\end{gather*}
By the definition of $\rho_{1}$
\begin{gather*}
2|\mathbf{E}(X_{j}+A_{j})A_{j-1}|\leq2\rho_{1}||X_{j}+A_{j}||_{2}%
||A_{j-1}||_{2}\\
\leq\rho_{1}^{2}\mathbf{E}(X_{j}+A_{j})^{2}+\mathbf{E(}A_{j-1})^{2}\text{ ,}%
\end{gather*}
which combined with the previous identity gives%
\[
\mathbf{E}[\mathbf{P}_{j}(S_{n}-S_{j-1})]^{2}\geq(1-\rho_{1}^{2}%
)\mathbf{E}(X_{j}+A_{j})^{2}\text{.}%
\]
Therefore, by summing these inequalities we obtain
\begin{equation}
\sigma_{n}^{2}\geq(1-\rho_{1}^{2})\sum_{i=1}^{n}\mathbf{E}(X_{i}+A_{i}%
)^{2}\text{.} \label{ineq1}%
\end{equation}
On the other hand, by (\ref{sigma}) and the properties of conditional
expectation
\[
\sigma_{n}^{2}=\sum_{i=1}^{n}\mathbf{E}(X_{i}+A_{i})^{2}-\sum_{i=1}%
^{n}\mathbf{E}A_{i-1}^{2}\text{ .}%
\]
By introducing this identity in relation (\ref{ineq1}) and changing the
variable of summation we obtain
\[
\sigma_{n}^{2}\geq(1-\rho_{1}^{2})[\sum_{i=1}^{n}\mathbf{E}A_{i}^{2}%
+\sigma_{n}^{2}]\text{ .}%
\]
Solving this inequality for $\sigma_{n}^{2}$ gives%
\begin{equation}
\sigma_{n}^{2}\geq\frac{1-\rho_{1}^{2}}{\rho_{1}^{2}}\sum_{i=1}^{n}%
\mathbf{E}A_{i}^{2}\text{ .} \label{est1}%
\end{equation}
(If $\rho_{1}=0,$ then $\sum_{i=1}^{n}\mathbf{E}A_{i}^{2}=0.)$

Starting now from $X_{i}=(X_{i}+A_{i})-A_{i}$, by the Cauchy-Schwarz
inequality, we have
\[
\mathbf{E}X_{i}^{2}\leq\mathbf{E}(X_{i}+A_{i})^{2}+\mathbf{E}A_{i}%
^{2}+2\mathbf{||}X_{i}+A_{i}||_{2}\mathbf{||}A_{i}||_{2}\text{ .}%
\]
We sum these inequalities, then we apply H\"{o}lder inequality and finally use
relations (\ref{ineq1}) and (\ref{est1}) and some simple calculations to
obtain
\begin{gather*}
b_{n}^{2}\leq\sum_{i=1}^{n}\mathbf{E(}X_{i}+A_{i})^{2}+\sum_{i=1}%
^{n}\mathbf{E}A_{i}^{2}+\\
2\left(  \sum_{i=1}^{n}\mathbf{E(}X_{i}+A_{i})^{2}\sum_{j=1}^{n}%
\mathbf{E}A_{j}^{2}\right)  ^{1/2}\leq\frac{1+\rho_{1}}{1-\rho_{1}}\sigma
_{n}^{2}\text{ .}%
\end{gather*}
Therefore
\[
\sigma_{n}^{2}\geq\frac{1-\rho_{1}}{1+\rho_{1}}b_{n}^{2}%
\]
and the lower bound is established.

We shall establish now the upper bound. By simple algebra, Cauchy-Schwarz and
H\"{o}lder inequalities, we have
\begin{align*}
\sigma_{n}^{2}  &  =-b_{n}^{2}+2\sum_{i=1}^{n}\mathbf{E[}X_{i}(X_{i}+A_{i})]\\
&  \leq-b_{n}^{2}+2b_{n}[\sum_{i=1}^{n}\mathbf{E}(X_{i}+A_{i})^{2}%
]^{1/2}\text{ .}%
\end{align*}
Since for any two positive numbers, $a$ and $b$, we have $2ab\leq(1-\rho
_{1})^{-1}a^{2}+(1-\rho_{1})b^{2}$, we obtain
\[
\sigma_{n}^{2}\leq\frac{\rho_{1}}{1-\rho_{1}}b_{n}^{2}+(1-\rho_{1})\sum
_{i=1}^{n}\mathbf{E}(X_{i}+A_{i})^{2}\text{ .}%
\]
This last inequality, combined with (\ref{ineq1})\ and solved for $\sigma
_{n}^{2}$ gives
\[
\sigma_{n}^{2}\leq\frac{1+\rho_{1}}{1-\rho_{1}}b_{n}^{2}\text{ ,}%
\]
and the upper bound is established. $\ \Diamond$

\begin{remark}
Notice that for an independent vector, $\rho_{1}=0$ and Proposition
\ref{lowerB2} can be viewed as an extension of the classical estimate for
variance in the independent case, $\mathbf{E}S_{n}^{2}=\sum_{i=1}%
^{n}\mathbf{E}X_{i}^{2}$.
\end{remark}

As a corollary we obtain the following result in terms of the coefficient of
contraction $\delta$ defined by (\ref{defdelta}) that improves the known
results in the literature (see for instance Section 1.2.2. in \cite{Io2} and
Proposition 3.2 in \cite{SV}).

\begin{corollary}
\bigskip\label{lowerB2delta}Let $(X_{1},X_{2},...,X_{n})$ be as in Proposition
\ref{lowerB2}. If $\delta_{1}<1$ then%
\[
\frac{1-\delta_{1}}{(1+\sqrt{\delta_{1}})^{2}}\sum_{i=1}^{n}\mathbf{E}%
X_{i}^{2}\leq\mathbf{E}S_{n}^{2}\leq\frac{(1+\sqrt{\delta_{1}})^{2}}%
{1-\delta_{1}}\sum_{i=1}^{n}\mathbf{E}X_{i}^{2}\text{ .}%
\]

\end{corollary}

\textbf{Proof}. It was established in Lemma 4.1 in \cite{SV} that $\rho
_{1,n}<\sqrt{\delta_{1,n}}.$

Then, since the function $(1-x)/(1+x)$ is decreasing, it follows by
Proposition \ref{lowerB2}~that
\[
\frac{1-\sqrt{\delta_{1}}}{1+\sqrt{\delta_{1}}}b_{n}^{2}\leq\mathbf{E}%
S_{n}^{2}\leq\frac{1+\sqrt{\delta_{1}}}{1-\sqrt{\delta_{1}}}b_{n}^{2}\text{ .}%
\]

$\Diamond$

\section{Proofs of the main results}

\subsection{Proof of Theorem \ref{CLTMarkov}}

We verify the conditions of Proposition \ref{markov}. Condition (\ref{C1})
follows easily by conditions (\ref{bound}) and (\ref{Dob}) combined with
Proposition \ref{lowerB2} in the following way:%
\begin{equation}
\frac{1}{\sigma_{n}^{4}}\sum_{j=1}^{n}\mathbf{E}X_{n,j}^{4}\leq\frac{C_{n}%
^{2}}{\sigma_{n}^{4}}\sum_{j=1}^{n}\mathbf{E}X_{n,j}^{2}\leq\frac{2C_{n}%
^{2}b_{n}^{2}}{\lambda_{n}b_{n}^{2}\sigma_{n}^{2}}=\frac{2C_{n}^{2}}%
{\lambda_{n}\sigma_{n}^{2}}\rightarrow0\text{ as }n\rightarrow\infty.
\label{four}%
\end{equation}
To verify (\ref{C2}) we fix $\varepsilon>0$ and start from%
\[
\mathbf{P}(\max_{1\leq j\leq n}|A_{j}|\geq\varepsilon\sigma_{n})\leq
\frac{\mathbf{E}\exp(t\max_{1\leq j\leq n}|A_{j}|)}{\exp t\varepsilon
\sigma_{n}}%
\]
For $t=\lambda_{n}/(6C_{n})$, by taking into account Lemma \ref{exp}, we
obtain \
\[
\mathbf{P}(\max_{1\leq j\leq n}|A_{j}|\geq\varepsilon\sigma_{n})\leq
(1+\frac{b_{n}}{3C_{n}})^{2}\exp(-\varepsilon\frac{\lambda_{n}\sigma_{n}%
}{6C_{n}})
\]
and then, (\ref{C2}) follows provided we verify
\begin{equation}
(1+\frac{b_{n}}{3C_{n}})\exp(-\varepsilon\frac{\lambda_{n}\sigma_{n}}{12C_{n}%
})\rightarrow0\text{ as }n\rightarrow\infty\text{ .} \label{conv}%
\end{equation}
Notice that assumption (\ref{Dob}) implies
\[
\frac{\lambda_{n}\sigma_{n}}{C_{n}}\rightarrow\infty\text{ as }n\rightarrow
\infty\text{ ,}%
\]
that further implies
\[
\exp(-\varepsilon\frac{\lambda_{n}\sigma_{n}}{12C_{n}})\rightarrow0\text{ as
}n\rightarrow\infty\text{ .}%
\]
Moreover%
\[
\frac{b_{n}}{C_{n}}\exp(-\varepsilon\frac{\lambda_{n}\sigma_{n}}{12C_{n}%
})=\exp\left(  \ln(\frac{b_{n}}{C_{n}})-\varepsilon\frac{\lambda_{n}\sigma
_{n}}{12C_{n}}\right)
\]
and (\ref{conv}) follows if we show that
\begin{equation}
\ln(\frac{b_{n}}{C_{n}})-\varepsilon\frac{\lambda_{n}\sigma_{n}}{12C_{n}%
}\rightarrow-\infty\text{ .} \label{inf}%
\end{equation}
We write now%
\[
\ln(\frac{b_{n}}{C_{n}})=\ln(\frac{\lambda_{n}\sigma_{n}}{C_{n}})+\ln
(\frac{b_{n}}{\lambda_{n}\sigma_{n}})\text{ .}%
\]
and notice that%
\[
\ln(\frac{\lambda_{n}\sigma_{n}}{C_{n}})-\varepsilon\frac{\lambda_{n}%
\sigma_{n}}{24C_{n}}\rightarrow-\infty\text{ ,}%
\]
so, in order for (\ref{inf}) to hold it is enough to show that for and $n$
sufficiently large%
\[
\ln(\frac{b_{n}}{\lambda_{n}\sigma_{n}})\leq\varepsilon\frac{\lambda_{n}%
\sigma_{n}}{24C_{n}}\text{ . }%
\]
This fact follows if%
\begin{equation}
\frac{C_{n}}{\lambda_{n}\sigma_{n}\ }\ln(\frac{b_{n}}{\lambda_{n}\sigma_{n}%
})\rightarrow0\text{ .} \label{convneeded}%
\end{equation}
Now, by Proposition \ref{lowerB2} we have%
\[
\frac{\lambda_{n}^{1/2}}{2}\leq\frac{b_{n}}{\sigma_{n}}\leq\frac{2}%
{\lambda_{n}^{1/2}}\text{.}%
\]
So, condition (\ref{convneeded}) is satisfied provided%
\[
\frac{\lambda_{n}\sigma_{n}}{C_{n}(1+|\ln(\lambda_{n})|)}\rightarrow
\infty\text{ .}%
\]
This is exactly the condition that we impose in (\ref{Dob}). Thus (\ref{C2}) holds.

We verify (\ref{C3}) by analyzing the variance of both terms involved. For the
first term we use Proposition \ref{lowerB2} together with conditions
(\ref{bound}) and (\ref{Dob}) and obtain
\begin{gather}
\frac{1}{\sigma_{n}^{4}}var\sum_{j=1}^{n}X_{n,j}^{2}\leq\frac{2}{\lambda
_{n}\sigma_{n}^{4}}\sum_{j=1}^{n}\mathbf{E}X_{n,j}^{4}\leq\label{varX}\\
\frac{2C_{n}^{2}b_{n}^{2}}{\lambda_{n}\sigma_{n}^{4}}\leq\frac{4C_{n}^{2}%
}{\lambda_{n}^{2}\sigma_{n}^{2}}\rightarrow0\text{ as }n\rightarrow
\infty.\nonumber
\end{gather}
To deal with the second term, first we apply Proposition \ref{lowerB2} to
estimate the variance, then we use condition (\ref{bound}), and finally we
take into account the inequality (\ref{est1}) and so,
\begin{gather*}
\frac{1}{\sigma_{n}^{4}}var\sum_{j=1}^{n}X_{n,j}A_{n,j}{\large \leq}\frac
{2}{\lambda_{n}\sigma_{n}^{4}}\sum_{j=1}^{n}\mathbf{E(}X_{n,j}^{2}A_{n,j}%
^{2})\\
\leq\frac{2C_{n}^{2}}{\lambda_{n}\sigma_{n}^{4}}\sum_{j=1}^{n}\mathbf{E}%
A_{n,j}^{2}\leq\frac{2C_{n}^{2}\sigma_{n}^{2}}{\lambda_{n}^{2}\sigma_{n}^{4}%
}=\frac{2C_{n}^{2}}{\lambda_{n}^{2}\sigma_{n}^{2}}%
\end{gather*}
which converges to $0$ under (\ref{Dob}). $\Diamond$

\subsection{Proof of Corollary \ref{Lindeberg}}

This corollary follows from Theorem \ref{CLTMarkov} via a truncation argument.

First construct $\varepsilon_{n}\rightarrow0$ slowly enough such that
condition (\ref{lindeberg}) is still satisfied. We truncate the variables at
the level $T_{n}=\varepsilon_{n}h(\lambda_{n})\sigma_{n},$ and denote%
\[
X_{n,i}^{^{\prime}}=X_{n,i}I(|X_{n,i}|\leq T_{n})-\mathbf{E}X_{n,i}%
I(|X_{n,i}|\leq T_{n})
\]
and
\[
X_{n,i}^{^{"}}=X_{n,i}-X_{n,i}^{^{\prime}}\text{ .}%
\]
We show that the contribution of $\sum_{i=1}^{n}X_{n,i}^{^{"}}/\sigma_{n}$ is
negligible in $\mathbf{L}_{2}$ and therefore is negligible for the convergence
in distribution. To estimate its variance we apply Proposition \ref{lowerB2}
and then we take into account the Lindeberg condition (\ref{lindeberg}). We
obtain
\begin{gather*}
\frac{1}{\sigma_{n}^{2}}var(\sum_{i=1}^{n}X_{n,i}^{^{"}})\leq\frac{2}%
{\lambda_{n}\sigma_{n}^{2}}\sum_{i=1}^{n}\mathbf{E(}X_{n,i}^{^{"}})^{2}\\
\leq\frac{4}{\lambda_{n}\sigma_{n}^{2}}\sum_{i=1}^{n}\mathbf{E}X_{n,i}%
^{2}I(|X_{n,i}|>T_{n})\rightarrow0\text{ as }n\rightarrow\infty\text{ .}%
\end{gather*}
Then, with the notation $(\sigma_{n}^{^{\prime}})^{2}=var$ $\sum_{i=1}%
^{n}X_{n,i}^{^{\prime}}$ we easily derive from the last convergence that
\[
\lim_{n\rightarrow\infty}\frac{1}{\sigma_{n}^{2}}(\sigma_{n}^{^{\prime}}%
)^{2}=1\text{ .}%
\]
Finally, we apply Theorem \ref{CLTMarkov} to $X_{n,i}^{^{\prime}}$ with
$\ C_{n}=2\varepsilon_{n}h(\lambda_{n})\sigma_{n}.$ We verify (\ref{Dob}) by
using the definition of $h(\lambda_{n})$, since%
\[
\frac{C_{n}(1+|\ln(\lambda_{n})|)}{\lambda_{n}\sigma_{n}^{^{\prime}}%
}=2\varepsilon_{n}\frac{\sigma_{n}}{\sigma_{n}^{^{\prime}}}\rightarrow0\text{
as }n\rightarrow\infty\text{ }.
\]
and the result follows. $\ \Diamond$

\subsection{Proof of Theorem \ref{Dobr}}

Using our tools we give a short proof of this theorem for completeness. We
verify conditions of Proposition \ref{markov}. Notice that under the
assumptions of this theorem Conditions (\ref{C1}) and (\ref{C3}) are verified
exactly as in the proof of Theorem \ref{CLTMarkov} by taking into account that
$\rho_{n,1}\leq\sqrt{\delta_{n,1}}$ and replacing Proposition \ref{lowerB2} by
its Corollary \ref{lowerB2delta}. The main difference is now that condition
(\ref{C2}) follows easily by the estimate
\[
||\mathbf{E}(X_{n,k}|\xi_{n,j})||_{\infty}\leq2\delta_{n,1}^{k-j}C_{n}\text{
a.s. }%
\]
This inequality implies
\[
\frac{1}{\sigma_{n}}||\mathbf{E}(S_{n}-S_{n,j}|\xi_{n,j})||_{\infty}\leq
\frac{1}{\sigma_{n}}\sum_{i=j+1}^{n}||\mathbf{E}(X_{n,i}|\xi_{n,j})||_{\infty
}\leq\frac{2C_{n}}{\alpha_{n}^{3/2}b_{n}}%
\]
which converges to $0$ as $n\rightarrow\infty$ \ by condition (\ref{CD}). The
theorem is established. $\Diamond$

\section{ Appendix}

In this section we estimate the moments and the exponential moments for the
quantity $A_{j}=A_{n,j}=\mathbf{E}(S_{n}-S_{j}|\mathcal{F}_{j})$ where $S_{j}$
are the partial sums associated to a vector of centered random variables
$(X_{j})_{1\leq j\leq n}$ defined on a probability space $(\Omega
,\mathcal{K},P)$, adapted to an increasing filtration of sub-sigma fields of
$\mathcal{K}$, $(\mathcal{F}_{j})_{1\leq j\leq n}$, $\mathcal{F}%
_{0}=\{\varnothing,\Omega\}$ and $\rho_{k}$ is defined by (\ref{defrogen}).

\begin{lemma}
\label{ppower}Let $p\geq2$ be a real number and assume the variables have
finite moments of order $p$. Then,%
\[
\sum_{j=1}^{n}\mathbf{E{\large |}}A_{n,j}|^{p}\leq2^{p-2}(\sum_{k=1}^{n-1}%
\rho_{k}^{2/p})^{p}\sum_{i=1}^{n}\mathbf{E}|X_{i}|^{p}\text{ .}%
\]
If for a certain $0<\rho<1$ we have $\rho_{k}\leq\rho^{k}$ then, for any
$p\geq2,$
\[
\sum_{j=1}^{n}\mathbf{E{\large |}}A_{n,j}|^{p}\leq p^{p}\frac{1}{4(1-\rho
)^{p}}\sum_{i=1}^{n}\mathbf{E}|X_{i}|^{p}\text{ .}%
\]

\end{lemma}

\textbf{Proof}. For simplicity, we shall drop the index $n$ from the notation.
For $j$ fixed, $1\leq j\leq n$, and $1\leq k\leq n-j$ let
\begin{equation}
a_{k}=a_{k}(j)=\frac{\rho_{k}^{2/p}}{\sum_{i=1}^{n-j}\rho_{i}^{2/p}}\text{ .}
\label{defa}%
\end{equation}
Notice that $\sum_{k=1}^{n-j}a_{k}=1$. By the fact that $x\rightarrow|x|^{p}$
is a convex function, we easily obtain
\begin{align*}
\mathbf{{\large |}}A_{j}|^{p}  &  =|\sum_{i=j+1}^{n}a_{i-j}a_{i-j}%
^{-1}\mathbf{E}(X_{i}|\mathcal{F}_{j})|^{p}\leq\sum_{i=j+1}^{n}a_{i-j}%
|a_{i-j}^{-1}\mathbf{E}(X_{i}|\mathcal{F}_{j})|^{p}\\
&  =\sum_{i=j+1}^{n}a_{i-j}^{1-p}|\mathbf{E}(X_{i}|\mathcal{F}_{j})|^{p}%
=\sum_{k=1}^{n-j}a_{k}^{1-p}(j)|\mathbf{E}(X_{j+k}|\mathcal{F}_{j})|^{p}\text{
.}%
\end{align*}
Then, since $a_{k}(j)\geq a_{k}(1)\ $and $1-p<0,$ it follows
\[
\mathbf{{\large |}}A_{j}|^{p}\leq\sum_{k=1}^{n-j}a_{k}^{1-p}(1)|\mathbf{E}%
(X_{j+k}|\mathcal{F}_{j})|^{p}\text{ .}%
\]
Next, we use the fact that by the interpolation theory (see Theorem 4.12 in
\cite{br}), for $p\geq2$,%
\[
\mathbf{E}|\mathbf{E}(X_{j+u}|\mathcal{F}_{j})|^{p}\leq2^{p-2}\rho_{u}%
^{2}\mathbf{E}|X_{j+u}|^{p}\text{ .}%
\]
Combining now these two facts and summing the relations, we obtain,
\begin{align*}
\sum_{j=1}^{n-1}\mathbf{E{\large |}}A_{j}|^{p}  &  \leq2^{p-2}\sum_{j=1}%
^{n-1}\sum_{k=1}^{n-j}a_{k}^{1-p}(1)\rho_{k}^{2}\mathbf{E}|X_{j+k}|^{p}\\
&  \leq2^{p-2}\sum_{k=1}^{n-1}a_{k}^{1-p}(1)\rho_{k}^{2}\sum_{j=1}%
^{n}\mathbf{E}|X_{j}|^{p}\text{.}%
\end{align*}
By (\ref{defa}) we notice that
\[
\sum_{k=1}^{n-1}a_{k}^{1-p}(1)\rho_{k}^{2}=(\sum_{k=1}^{n-1}\rho_{k}%
^{2/p})^{p}%
\]
and the first part of this lemma follows.

For proving the second part of this lemma we take into account that a simple
computation based on the fact that $1-x^{\beta}\geq\beta(1-x)$ for $0<x\leq1$
and $0<\beta\leq1$ gives%
\[
\sum_{k=1}^{n-1}\rho_{k}^{2/p}=\sum_{k=1}^{n-1}\rho^{2k/p}\leq\frac{1}%
{1-\rho^{2/p}}\leq\frac{p}{2}\frac{1}{(1-\rho)}\text{\ ,}%
\]
which combined with the first part of the lemma gives the result. $\Diamond$

For the next lemma we recall the definition (\ref{def b}).

\begin{lemma}
\label{exp}Assume that there is $C>0$ such that $\max_{1\leq k\leq n}%
|X_{k}|\leq C$ $\ a.s.$ and for a certain $0<\rho<1$ we have $\rho_{k}\leq
\rho^{k}.$ Then, for any nonnegative $t\leq(1-\rho)/(6C)$%
\[
\mathbf{E}\exp(t\max_{1\leq j\leq n}|A_{j}|)\leq(1+\frac{2tb_{n}}{1-\rho}%
)^{2}\text{ .}%
\]
In particular for $t=\frac{1-\rho}{6C}$ we have
\[
\mathbf{E}\exp(\frac{1-\rho}{6C}\max_{1\leq j\leq n}|A_{j}|)\leq(1+\frac
{b_{n}}{3C})^{2}\text{ .}%
\]

\end{lemma}

\textbf{Proof}. We start the estimate by the Taylor expansion and majorate the
maximum term by the sum:%
\begin{gather*}
\ \mathbf{E}\exp(t\max_{1\leq j\leq n}|A_{j}|)\leq1+\sum_{p=2}^{\infty}%
\frac{t^{p}}{p!}\mathbf{E}\max_{1\leq j\leq n}|A_{j}|^{p}+t\mathbf{E}%
\max_{1\leq j\leq n}|A_{j}|\leq\\
1+\sum_{p=2}^{\infty}\sum_{j=1}^{n}\frac{t^{p}}{p!}\mathbf{E|}A_{j}%
|^{p}+t\mathbf{E}\max_{1\leq j\leq n}|A_{j}|=I+II\text{ ,}%
\end{gather*}
where
\[
I=1+\sum_{p=2}^{\infty}\sum_{j=1}^{n}\frac{t^{p}}{p!}\mathbf{E|}A_{j}%
|^{p}\text{ .}%
\]
By Lemma \ref{ppower} and because by the Stirling approximation we have
$p^{p}\leq3^{p-1}p!$ for $p\geq2$, we obtain%
\[
\sum_{j=1}^{n}\mathbf{E{\large |}}A_{j}|^{p}\leq\frac{p^{p}}{4(1-\rho)^{p}%
}\sum_{i=1}^{n}\mathbf{E}|X_{i}|^{p}\leq\frac{3^{p}p!C^{p-2}b_{n}^{2}%
}{12(1-\rho)^{p}}\text{ .}%
\]
Introducing this estimate in the expression of $I$ we have
\[
I\leq1+\sum_{p=2}^{\infty}\frac{(3t)^{p}C^{p-2}b_{n}^{2}}{12(1-\rho)^{p}%
}=1+\frac{9t^{2}b_{n}^{2}}{12(1-\rho)^{2}}\sum_{p=2}^{\infty}\frac
{(3tC)^{p-2}}{(1-\rho)^{p-2}}\text{ .}%
\]
For $t\leq\frac{1-\rho}{6C}$ we easily derive
\[
I\leq1+\frac{3t^{2}b_{n}^{2}}{4(1-\rho)^{2}}\left(  1-\frac{3tC}{1-\rho
}\right)  ^{-1}\leq1+\frac{3t^{2}b_{n}^{2}}{2(1-\rho)^{2}}\text{.}%
\]
Moreover, by Lemma \ref{ppower} it follows that
\[
II=t(\mathbf{E}\max_{1\leq j\leq n}|A_{j}|)\leq t\left(  \sum_{j=1}%
^{n}\mathbf{E}A_{j}^{2}\right)  ^{1/2}\leq\frac{tb_{n}}{1-\rho}%
\]
and overall%
\[
I+II\leq1+\frac{3t^{2}b_{n}^{2}}{2(1-\rho)^{2}}+\frac{tb_{n}}{1-\rho}%
\leq(1+\frac{2tb_{n}}{1-\rho})^{2}%
\]
and the lemma is established. \ $\Diamond$

\section{Acknowledgement}

The author is grateful to Richard Bradley, Sergey Utev and Sunder Sethuraman
for useful discussions on the subject. Many thanks go to the referee for
carefully reading the manuscript and for numerous suggestions that improved
the presentation of this paper.

\end{document}